\documentclass[12pt]{amsart}

\usepackage{amssymb}

\newtheorem{theorem}{Theorem}[section]

\newtheorem{lemma}[theorem]{Lemma}
\newtheorem{corollary}[theorem]{Corollary}
\newtheorem{defn}[theorem]{Definition}
\newtheorem{remark}[theorem]{Remark}

\numberwithin{equation}{section}

\numberwithin{equation}{section}

\begin{document}
\baselineskip=16.5pt

\title[On the Spectrum of Asymptotic Slopes]{On the Spectrum of
Asymptotic Slopes}

\author{A. J. Parameswaran}

\address{School of Mathematics, Tata Institute of Fundamental
Research,  Homi Bhabha Road, Mumbai 400005, India}

\email{param@math.tifr.res.in}

\email{subramnn@math.tifr.res.in}

\author{S. Subramanian}

\maketitle

\begin{abstract}
The slopes of maximal subbundles of rank $s$ divided by the degree of the
map under various pull backs form a bounded collection of numbers called
the $s$-spectrum of the bundle. We study the supremum of the $s$-spectrum
and determine it in terms of the Harder Narasimhan filtration of the
bundle. 
\end{abstract}

\section{Introduction}

Line subbundles of maximal degree in a rank $2$ bundle on a curve have
been studied in [LN] and many other subsequent papers. In [BP] it was
shown that a vector bundle is strongly semistable if and only if the slope
of the
maximal subbundle of a given rank in a finite pull back is bounded by
the slope of the pull back bundle.

In this note we study the behaviour of maximal subbundles of vector
bundles on curves after finite pull backs. These slopes of maximal
subbundles of rank $s$ divided by the degree of the map under various
pull backs form a bounded collection of numbers called the
$s$-spectrum of the bundle. We study the supremum of the $s$-spectrum
and determine it in terms of the Harder Narasimhan filtration of the
bundle. We also give a criterion for a spectrum value to be isolated.



\section{Preliminaries} 

Let $k$ be an algebraically closed field of arbitrary characteristic.
Let $C$ be an irreducible smooth projective curve over $k$ of genus
$g$. For a vector bundle $V$ over $C$ we will denote by ${\mathbb
G}r(s,V)$ the Grassmann bundle over $C$ defined by the space of all
$s$ dimensional quotients of $V$. The tautological line bundle over
${\mathbb G}r(s,V)$ will be denoted by ${\mathcal O}_{{\mathbb
G}r(s,V)}(1)$. This is the determinant of the universal quotient
bundle on ${\mathbb G}r(s,V)$. Tensor powers of this tautological line
bundle are denoted by ${\mathcal O}_{{\mathbb G}r(s,V)} (n)$. 
The slope of $V$ is defined
as $$\mu(V):= \frac{{\rm deg}(V)}{{\rm rank}(V)}$$ $V$ is called {\em
semistable} if for all subbundles $W\subset V$ of positive rank
$\mu(W)\leq \mu(V)$.  $V$ is called {\em strongly semistable} if the
iterated Frobenius pull back, ${F_C^n}^*(V)$, is semistable for all
$n>0$, where $F_C:C\to C$ denotes the Frobenius morphism. If $V$ is a
strongly semistable vector bundle of rank $r$, then for any
representation of $GL_r$ into $GL_n$, the induced rank $n$ vector
bundle is strongly semistable (cf. [RR]).

Given a vector bundle $V$, there is a unique filtration (the Harder
Narasimhan filtration) $V_{\bullet}:=\{0=V_0\subset \cdots \subset
V_l=V\}$ such that each $V_i/V_{i-1}$ is semistable and the slopes of
successive quotients are strictly decreasing, i.e., $\mu_i:=\mu(V_i/
V_{i-1})>\mu_{i+1}:=\mu(V_{i+1}/V_i)$. The subbundles $V_i$ are
defined inductively as the inverse image of the maximal subbundle of
maximal slope in $V/V_{i-1}$. We call the bundles $V_i/V_{i-1}$, the
{\em Harder Narasimhan factors} of the bundle $V$. The bundle $V_1$ is
called the {\em maximal subbundle} of $V$ and denoted by $V_{\rm max}$.
It's slope $\mu(V_1)=\mu(V_{\rm max})$ is called the maximal slope of
$V$ and denoted by $\mu_{\rm max}(V)$. We recall the following
result whose proof we will omit. 

\begin{lemma} If $V$ and $W$ are semistable vector bundles over a
smooth curve $C$, with $\mu(V)>\mu(W)$, then $Hom_C(V,W)= H^0(C,
V^*\otimes W)=0$.
\end{lemma} 

Now  we have  the   following  result  which  gives  enough  complete
intersection curves in the projective bundle.

\begin{lemma}\label{lemma1} Let $V$ be a strongly semistable vector
bundle of rank $r$ over $C$ and $L$ a line bundle of degree $>
2g-ns\mu(V)$ with $n>0$. Then ${\mathcal O}_{{\mathbb G}r
(s,V)}(n)\otimes \pi^*L$ separates points on ${\mathbb G}r(s,V)$.
\end{lemma}

\begin{proof} Since $\pi:{\mathbb G}r(s,V)\to C$ is a smooth fibration
with Grassmann varieties as fibres and higher cohomologies of ample
bundles vanish on Grassmannians, it follows that $$H^1({\mathbb
G}r(s,V)~,{\mathcal O}_{{\mathbb G}r (s,V)}(n)\otimes \pi^*L) \cong
H^1 (C, \pi_*({\mathcal O}_{{\mathbb G}r(s,V)}(n)\otimes \pi^*L))\cong
H^1 (C, V_{s,n} \otimes L) $$ where $V_{s,n}$ is the vector bundle
associated to $V$ by the Weyl module with highest weight $n\omega_s$
with $\omega_s$ as the fundamental weight corresponding to the $s$-th
exterior power representation. Notice that the slope of $\wedge^s V$
is equal to $s\mu(V)$.

Hence for any two points $x,y\in C$, Serre duality implies
$$H^1(C,~V_{s,n}\otimes L(-x-y))~=~H^0(C,~V_{s,n}^*\otimes L^{-1}(x+y)
\otimes K_C)^*$$ To show this cohomology group vanishes, it suffices
to show the semistable bundle $V_{s,n}^*\otimes L^{-1}(x+y) \otimes
K_C$ has negative slope. Now we have: $$\mu(V_{s,n}^*\otimes
L^{-1}(x+y)\otimes K_C) ={\mu}(V_{s,n}^*) - {\rm deg}~L + 2g = -ns\mu
(V) - {\rm deg}~L +2g $$ which is negative if and only if ${\rm
deg}~L> 2g-ns{\mu}(V)$.  Hence for two distinct points $x,y\in C$, $$
H^0 (C,~V_{s,n}\otimes L) \to H^0(C,~V_{s,n}\mid_{x+y})$$ is
surjective whenever deg $ L> 2g-ns{\mu}(V)$. This implies ${\mathcal
O}_{{\mathbb G}r(s,V)}(n)\otimes \pi^*L$ separate points and surjects
onto sections of ${\mathcal O}_{{\mathbb G}r(V_x)}(n)$ on the
Grassmannian ${\mathbb G}r(s,V_x)$.  \end{proof}

Consider the universal exact sequence on ${\mathbb G}r(s,V)$:
\begin{equation}\label{opV1}
0\to{\mathcal S}\to\pi^*V\to {\mathcal Q}\to 0
\end{equation}

Then ${\rm det}~{\mathcal Q}\cong {\mathcal O}_{{\mathbb G}r(s,V)}(1)$
is the tautological line bundle. For any line bundle on ${\mathbb
G}r(s,V)$ let [~-~] denote the corresponding cycle class.  The Chow
group of $0$ cycles on ${\mathbb G}r(s,V)$ is canonically isomorphic to
the Chow group of $0$ cycles of the curve $C$ (cf. [F], Prop. 14.6.5,).
Then we have the following:

\begin{lemma}\label{chowrelations} Given any line bundle $L$ on $C$ and
any fibre $F$ of $\pi:{\mathbb G}r(s,V)\to C$ we have

\begin{equation}\label{O1tothetop} 
[{\mathcal O}_{{\mathbb G}r(s,V)}(1)]^{s(r-s)+1} = (s(r-s)+1)s\mu(V)
([{\mathcal O}_{{\mathbb G}r(s,V)}(1)]^{s(r-s)}\cdot F)
\end{equation}
\begin{equation}\label{O1tothetoponfibre} 
[{\mathcal O}_{{\mathbb G}r(s,V)}(1)]^{s(r-s)} \cdot [\pi^*L] = {\rm
deg L}([{\mathcal O}_{{\mathbb G}r(s,V)}(1)]^{s(r-s)}\cdot F) ~{\rm
and}~[\pi^*L]\cdot [\pi^*L]=0
\end{equation}
\end{lemma}

\begin{proof} If we pull back the bundle under a finite map, both sides of
the formulae gets multiplied by the degree of the map. Hence we may assume,
after a finite pull back if necessary,  that  there exists a line bundle
$L$ such that $L^{\otimes r}\cong {\rm det} V$., i.e., det $V$ has an
$r^{\rm th}$ root. Then $\oplus_1^r L\cong {\mathcal O}_C^{\oplus
r}\otimes L$. Let ${\mathcal L}^{-1}$ be a very ample line bundle on $C$
such that $V\otimes {\mathcal L}^{-1}$ and $L\otimes {\mathcal L}^{-1}$ 
are globally generated and have vanishing first cohomology.  

Let $Quot_P({\mathcal O}_C^{\oplus N}\otimes {\mathcal L})$ denote the
quot scheme over $C$ of quotient of fixed Hilbert polynomial $P$ (the
degree and rank determine the polynomial) of  the trivial vector bundle of
rank $N>> 0$ twisted by the line bundle ${\mathcal  L}$. Let $U\subset
Quot_P({\mathcal O}_C^{\oplus N}\otimes {\mathcal L})$ be the open set
where the universal quotient sheaf is locally free and $H^0(C,V\otimes
{\mathcal L}^{-1})\cong H^0(C,{\mathcal O}_C^{\oplus N})$. Then by [Ne]
(Remark 5.5, page 140) $U$ is smooth and irreducible and hence connected.
Then $\oplus_1^r L\cong {\mathcal O}_C^{\oplus r}\otimes L$ and $V$ belong
to the same quot scheme, in fact $U$.

Let ${\mathcal V}$ be the universal bundle over $C\times U$ and ${\mathbb
G}r(s,{\mathcal V})$ by the corresponding Grassmannian bundle over
$C\times U$. For each vector bundle $W$ representing a closed point
$[W]\in U$, the restriction ${\mathcal V}\mid_{C\times [W]}$ is
canonically isomorphic to  $W$. Hence the restriction of the Grassmannian
bundle ${\mathbb G}r(s,{\mathcal V})$ is canonically ${\mathbb G}r(s,W)$.
Since the isomorphism of $CH_0({\mathbb G}r(s,V))$ with $CH_0(C)$ is
canonical, it suffices to check the formula for any closed point (vector
bundle) of the open subset $U$ of the quot scheme $Quot_P({\mathcal
O}_C^{\oplus N}\otimes {\mathcal L})$. 

Hence we obtain,

$[{\mathcal O}_{{\mathbb G}r(s,V)}(1)]^{s(r-s)+1} 
\cong 
[{\mathcal O}_{{\mathbb G}r(s,\oplus{L})}(1)]^{s(r-s)+1} 
\cong
[{\mathcal O}_{{\mathbb G}r(s,\oplus{\mathcal O}_C\otimes
L)}(1)]^{s(r-s)+1} \\
\cong
\{[{\mathcal O}_{{\mathbb G}r(s,\oplus{\mathcal O}_C)}(1)]+ [L]^{\otimes
s}\}^{s(r-s)+1} 
= 
(s(r-s)+1)s {\rm deg} L ({\mathcal O}_{{\mathbb G}r(s,\oplus{\mathcal
O}_C)}(1).F)  \\
= (s(r-s)+1)s \mu(V) ({\mathcal O}_{{\mathbb G}r(s,\oplus L)}(1).F) 
= (s(r-s)+1)s \mu(V) ({\mathcal O}_{{\mathbb G}r(s,V)}(1).F)
$
\end{proof}

Choose constants $0<\epsilon_n\leq 1$ and line bundles $L_n$ on $C$
such that ${\rm deg}~L_n=2g-ns\mu(V)+\epsilon_n$. Then ${\mathcal
O}_{{\mathbb G}r(s,V)}(n)\otimes \pi^*L_n$ separates points by
Lemma~\ref{lemma1} and hence defines enough smooth complete
intersection curves in ${\mathbb G}r(s,V)$ by Bertini's
theorem.  Now we have the following result.

\begin{lemma}\label{lemma2} Let $D$ be an irreducible complete
intersection curve in ${\mathbb G}r(s,V)$ defined by sections of
${\mathcal O}_{{\mathbb G}r(s,V)}(n)\otimes \pi^* L_n$. Then 
$$
{\mu}({\mathcal S}\mid_D) = n^{s(r-s)}([{\mathcal O}_{{\mathbb
G}r(s,V)}(1)]^{s(r-s)}\cdot F) (\mu(V)-\frac{s(2g+\epsilon_n)}{n}) $$
\end{lemma}

\begin{proof} Since $D$ is a complete intersection, the degree of
${\mathcal O}_{{\mathbb G}r(s,V)}(1)$ on $D$ can be calculated as the
cup product (denoted by $\cdot$ in the Chow ring) of the cycle classes
of the corresponding divisors (line bundles) with the class of
${\mathcal O}_{{\mathbb G}r(s,V)}(1)$. Note that $[{\mathcal
O}_{{\mathbb G}r(s,V)}(n)] = n[{\mathcal O}_{{\mathbb G}r(s,V)}(1)]$
as the tensor product of line bundles gives the sum of the
corresponding classes.
Hence we can interpret ${\rm deg}~{\mathcal O}_{{\mathbb
G}r(s,V)}(1)\mid_D$ as
$${\rm deg}~{\mathcal O}_{{\mathbb
G}r(s,V)}(1)\mid_D ~=~
[D]\cdot [{\mathcal O}_{{\mathbb G}r(s,V)}(1)]~=~
([{\mathcal O}_{{\mathbb G}r(s,V)}(n)]+[\pi^*L_n])^{s(r-s)}\cdot
[{\mathcal
O}_{{\mathbb G}r(s,V)}(1)] $$
$$ = \{[{\mathcal O}_{{\mathbb G}r(s,V)}(n)]^{s(r-s)}+
s(r-s)[{\mathcal O}_{{\mathbb G}r(s,V)}(n)]^{s(r-s)-1}\cdot
[\pi^*L_n]\}\cdot
[{\mathcal O}_{{\mathbb G}r(s,V)}(1)] $$ 
$$ = n^{s(r-s)} [{\mathcal
O}_{{\mathbb G}r(s,V)}(1)]^{s(r-s)+1} + s(r-s)n^{s(r-s)-1}
[{\mathcal O}_{{\mathbb
G}r(s,V)}(1)]^{s(r-s)}\cdot [\pi^*L_n] 
$$ 

Now the degree of the universal subbundle on $D$,  
$${\rm deg}~{\mathcal S}\mid_D={\rm deg}~{\pi^*V}\mid_D~-
~{\rm deg}~{\mathcal Q}\mid_D$$
$$
= n^{s(r-s)}([{\mathcal O}_{{\mathbb G}r(s,V)}(1)]^{s(r-s)}\cdot F)~
{\rm deg}~V~- ~n^{s(r-s)}(s(r-s)+1)s\mu(V)
([{\mathcal O}_{{\mathbb G}r(s,V)}(1)]^{s(r-s)}\cdot F)$$
$$ - 
n^{s(r-s)-1}s(r-s){\rm deg}~L_n~([{\mathcal O}_{{\mathbb
G}r(s,V)}(1)]^{s(r-s)}\cdot F) $$

Now  by substituting for deg $L_n$ in
the above and simplifying this expression we get

$$\mu({\mathcal S}\mid_D) = n^{s(r-s)}([{\mathcal O}_{{\mathbb
G}r(s,V)}(1)]^{s(r-s)}\cdot F)
(\mu(V)-\frac{s(2g+\epsilon_n)}{n}) $$

\end{proof}

Note that the degree of $\pi_D:D\to C$ is equal to cardinality of a general
fibre of $\pi_D$ over $x$ which equals $[{\mathcal O}_{{\mathbb
G}r(s,V_x)}(n)]^{s(r-s)}\cdot F=n^{s(r-s)}([{\mathcal O}_{{\mathbb
G}r(s,V)}(1)]^{s(r-s)}\cdot F)$. 

\section{Genuinely ramified maps}

Let us begin with the following definition.

\begin{defn} Let $f: D\to C$ be a finite morphism of integral
curves. Then $f$ is said to be {\em genuinely ramified} if $f$ is
separable and does not factor through an \'etale cover of $C$. 
\end{defn}

\begin{lemma}\label{genram}
A separable morphism $f:D\to C$ is genuinely ramified if and only if
$(f_*{\mathcal O}_{D})_{\rm max}\cong {\mathcal O}_{C}$.  
\end{lemma}

\begin{proof} By projection formula we have $$H^0(D, f^*W) = H^0(C,
W\otimes f_*{\mathcal O}_{D}) $$ Hence for a semistable bundle $S$ on
$C$ of positive slope $$Hom (S, f_*{\mathcal O}_{D}) = Hom (f^*S,
{\mathcal O}_{D}) =0 $$  because $f^*S$ remains semistable
of positive slope as $f$ is separable. This shows that $\mu_{\rm max}
f_*{\mathcal O}_{D}=0$.

If rank $(f_*{{\mathcal O}_{D}})_{\rm max} >1$, then it forms a sheaf of
subalgebras of $f_*{\mathcal O}_{D}$ on $C$. Hence by taking the
spectrum of $f_*{{\mathcal O}_{D}}_{\rm max}$ we obtain an \'etale cover
of $C$ factoring $f:D\to C$. Hence $f$ is not genuinely ramified. 
\end{proof}

In [N], Madhav Nori has constructed the fundamental group scheme
$\pi(X)$ of any complete variety $X$. This is constructed as a Tannaka
category whose objects are essentially finite vector bundles on $X$.
Further he has also constructed a principal $\pi(X)$ bundle
$\tilde{X}\to X$.  Let $\pi_1^{\rm alg}(X)$ denote the \'etale
fundamental group of a complete scheme $X$. Then one can show that
$\pi_1^{\rm alg}(X)$ is a quotient of $\pi(X)$ whose objects in the
Tannaka category are the vector bundles that are trivial on a finite
\'etale cover of $X$. This quotient morphism induces a $\pi_1^{\rm
alg}(X)$-bundle ${\mathcal X}\to X$ over $X$.

\begin{lemma}\label{genraminpi1} A separable morphism of curves $D\to C$
is
genuinely ramified if and only if the induced map $\pi_1^{\rm
alg}(D)\to \pi_1^{\rm alg}(C)$ is an epimorphism.  \end{lemma}

\begin{proof} Clearly an epimorphism on the fundamental group implies
genuine ramification. This follows from the fact that for any finite
\'etale morphism $f:D\to C$ of degree $d$ the index of the image of
the \'etale fundamental group is equal to $d$.

Any finite morphism $f:D\to C$ induces the map on \'etale fundamental
groups whose image $f_*(\pi_1^{\rm
alg}(D))=\Gamma$ is of finite index. Define $\tilde{D}:= {\mathcal
C}/\Gamma $. Then the induced map $\tilde{D}\to C$ is \'etale such
that $f:D\to C$ factors through this. This proves the converse.
\end{proof}

\begin{lemma}\label{genramofCIcurves} If $D$ is a general complete
intersection curve in ${\mathbb G}r(s,V)$ defined by an ample line
bundle, then the induced projection $D\to C$ is genuinely ramified
\end{lemma}

\begin{proof} Notice that the algebraic fundamental group of ${\mathbb
G}r(s,V)$ is naturally isomorphic to the algebraic fundamental group
of $C$ as Grassmannians are algebraically simply connected. If $D$ is a
complete intersection curve in ${\mathbb G}r(s,V)$, then the algebraic
fundamental group of $D$ surject onto the algebraic fundamental group
of ${\mathbb G}r(s,V)$ by the algebraic analogue of Lefshetz
Theorem. Hence $\pi_1^{\rm alg}(D)\to \pi_1^{\rm alg}(C)$ is
surjective. Now the result follows from Lemma~\ref{genraminpi1}.
\end{proof}

\begin{lemma}\label{stabletostable} Let $f:D\to C$ be genuinely
ramified morphism of smooth projective curves. Then\\ (a) If $V$ and
$W$ are two semistable bundles on $C$ of same slope, then
$$Hom_C(V,W)\cong Hom_D(f^*V,f^*W)$$ (b) If $V$ is a stable bundle on
$C$, then $f^*V$ is stable on $D$.\\ (c) If $V$ is a semistable bundle
on $C$ and $F\subset f^*V$ is a subbundle of same slope as $f^*V$,
then $F$ is isomorphic to the pull back of a subbundle of $V$. \\
\end{lemma}

\begin{proof} (a) Given two semistable bundles $V,~W$ of same slope on
$C$, we have $$Hom_D(f^*V,f^*W) \cong Hom_C(V, f_*f^*W)\cong Hom_C(V,
W\otimes f_*{\mathcal O}_{D}) \cong Hom_C(V,W)$$ The last equality
follows from the fact that $Hom_C(F, f_*{\mathcal O}_{D}/{\mathcal
O}_{C})=0$ for any semistable bundle $F$ of slope $\geq 0$ as genuine
ramification of $f$ implies that $f_*{\mathcal O}_{D}/{\mathcal
O}_{C}$ has negative maximal slope (see Lemma~\ref{genram}).

(b) Since the socle (maximal subbundle that is a direct sum stable
bundles (cf. [MR])) is unique, it follows that the socle of $f^*V$
descends to the socle of $V$ when $f$ is separable. Since $V$ is
stable, this descended bundle has to be $V$ itself. This shows that
the pull back of a stable bundle is polystable under any finite
separable map.  Now the stability of $f^*V$ for genuinely ramified
maps follows from (a) since projections to direct summands are
endomorphisms which do not come from below.

(c) Let $V$ be a semistable bundle over $C$. Let $F\subset f^*V$ be a
subbundle of same slope. Then the socle $S_F$ of $F$ is contained in
the socle $S_{f^*V}$ of $f^*V$ and hence a direct summand of
$S_{f^*V}$. But by uniqueness of the socle, $S_{f^*V}$ is $f^*(S_V)$.
Since stable bundles pull back to stable bundles, $S_F$ coincides with
some factors of $f^*(S_V)$ and hence is a pull back. Now the assertion
follows by induction on the rank applied to the bundle $F/{S_F}
\subset {f^*V}/{S_F}$.  
\end{proof}

\begin{corollary}\label{rank} Let $f:D\to C$ be a finite separable 
morphism of smooth projective irreducible curves. 
Then for any semistable vector bundle $W$ on $D$ we have\\ 
(a) $\mu_{\rm max}(f_*W)\leq \frac{\mu(W)}{{\rm deg}~f} $ \\
(b) If $\mu_{\rm max}(f_*W)~=~\frac{\mu(W)}{{\rm deg}~f}$, then
 ${\rm rank} (f_*W)_{\rm max} \leq {\rm rank}~W
\cdot {\rm rank}(f_*{\mathcal O}_D)_{\rm max} $ 
\end{corollary}

\begin{proof} First assertion follows from the fact that
$Hom_C(F,f_*W) \cong Hom_D(f^*F, W)$. Hence semistable bundles of
slope $>\frac{\mu(W)}{{\rm deg}~f}$ have no morphism to $f_*W$. 

If ${\rm rank}(f_*{\mathcal O}_D)_{\rm max}>1$, then the morphism
$f:D\to C$ factors through an \'etale morphism $\pi:\tilde C\to C$
such that $\tilde f: D\to \tilde C$ satisfies ${\tilde f}_*{\mathcal
O}_D ={\mathcal O}_{\tilde C}$ (see Lemma~\ref{genram}).

If $W'$ is any semistable bundle on $\tilde C$, then $\chi (W')= \chi
(\pi_*W')$ as the cohomologies do not change by taking direct images
under finite maps. By Riemann-Roch theorem we get $$ \chi (W') = ({\rm
rank}~(W'))(1-g_{\tilde C}+\mu(W'))= ({\rm rank}~(W'))({\rm
deg}~(\pi)(1-g_C)+\mu(W'))$$ and $$ \chi (\pi_*(W')) = {\rm
rank}~(\pi_*(W'))(1-g_C+\mu(\pi_*(W'))) = ({\rm deg}~(\pi){\rm
rank}~(W'))(1-g_C+\mu(\pi_*(W')))$$ since $\pi$ is \'etale,
$1-g_{\tilde C} = ({\rm deg}~\pi)(1-g_C)$.  Hence by comparing the
terms above we conclude that $\mu(W')=({\rm deg}~(\pi))
(\mu(\pi_*(W')))$. That $\pi_*W'$ is semistable of slope
$\frac{\mu(W')}{{\rm deg}~\pi}$ follows from (a).

Hence it suffices to prove (b) for the case when $f$ is genuinely ramified.

On the contrary, assume there is a semistable bundle $W$ with $\mu_{\rm
max} (f_*W) = {\frac{\mu(W)}{{\rm deg}~f}}$ and ${\rm rank}
(f_*W)_{\rm max} > {\rm rank} W$, then consider
the natural map, $f^*(f_*W)_{\rm max}\to Q\subset W$ with  image $Q$.
By taking direct
image we obtain $(f_*Q)_{\rm max} = (f_*W)_{\rm max}$. Since
$f^*(f_*W)_{\rm max}$ and $W$ are semistable of same slope, $Q$ is a
semistable vector bundle of same slope ${\frac{\mu(W)}{{\rm
deg}~f}}$. By Lemma~\ref{stabletostable}(c), $Q\cong f^*Q'$
is itself a pull back. Then $(f_*Q)_{\rm max} = (f_*f^*Q')_{\rm max} =
Q'$, hence has the same rank as $Q$, which is at most ${\rm rank}~W$,
a contradiction.  \end{proof}

\section{Asymptotic Slopes and Strong Semistability}

Let $C$ be a smooth curve defined over an algebraically closed field
$k$ of arbitrary characteristic.  Let $V$ be a vector bundle of rank
$r$ over $C$. For each $1\leq s< r$, we denote the slope of maximal
subbundle by $e_s(V)$.  $$ e_s(V) := Max~ \{~ \frac{{\rm deg} (W)}{s}
\mid W\subset V ~{\rm~is~a~subbundle~of~rank}~s~ \} $$ Define the
asymptotic $s$-spectrum ${\mathcal {AS}}_s(V)$ and the asymptotic
$s$-slope $\nu_s(V)$ as follows: $${\mathcal {AS}}_s(V) :=
\{~~\frac{e_s(f^*(V))}{{\rm deg}~f}~~\} $$
$$\nu_s(V):= Limsup~~~~\frac{e_s(f^*(V))}{{\rm deg}~f }
= Limsup ~~~~{\mathcal {AS}}_s(V) $$ where the supremum is taken over
all finite morphisms $f:D\to C$. Now we have the following criterion
for strong semistability in terms of the asymptotic slopes.

\begin{theorem}\label{theorem1} 
A vector bundle $V$ is strongly semistable if and only
if $\nu_s(V) = \mu(V) $ for some $s$. Then the asymptotic $s$-slopes
$\nu_s(V)$ are equal to the usual slope $\mu(V)$ for all $s$.
\end{theorem}

\begin{proof} Let $s$ be a given integer such that $1\leq s< r$.  In
[BP], it is proved that a vector bundle is strongly semistable if and only if
for every morphism $f:D\to C$, and for every subbundle $W\subset f^*V$
of rank $s$, the slope of $W$ is at most the slope of $f^*V$. Hence if
$V$ is not strongly semistable then there exists a finite morphism
$f:D\to C$ and a subbundle $W\subset f^*V$ of rank $s$ such that
$\mu(W)>\mu(f^*(V))$. Hence $\nu_s(V)>\mu(V)$.

Assume $V$ is strongly semistable. Then for any given finite map
$f:D\to C$, $f^*V$ is semistable and hence for every subbundle
$W\subset f^*V$ of rank $s$, the slope of $W$ is at most the slope of
$f^*V$.  Hence $\nu_s(V)\leq \mu(V)$.

Given $\epsilon >0$, choose $n\gg 0$ such that
$s(2g+\epsilon_n)<n\epsilon$. Then for the line bundle $L_n$,
${\mathcal O}_{{\mathbb G}r(s,V)}(n)\otimes \pi^*L$ separate points by
Lemma~\ref{lemma1}. The kernel of the universal sequence on a general
complete intersection $D$ on ${\mathbb G}r(s,V)$ defined by ${\mathcal
O}_{{\mathbb G}r(s,V)}(n)\otimes \pi^*L$ (which exists as it separates
points) has slope $n^{s(r-s)}([{\mathcal O}_{{\mathbb
G}r(s,V)}(1)]^{s(r-s)}\cdot F) (\mu(V)-\frac{s(2g+\epsilon_n)}{n})$,
by Lemma~\ref{lemma2}.  The degree of $\pi_D:D\to C$ is
$n^{s(r-s)}([{\mathcal O}_{{\mathbb G}r(s,V)}(1)]^{s(r-s)}\cdot F)$.
Hence dividing by the degree of $D\to C$, we obtain a number whose
difference with $\mu(V)$ is less than $\epsilon$.  Hence
$\nu_{r-s}(V)=\mu(V)$ for semistable bundles.
\end{proof}

\begin{remark}\label{limitslope} Let $V$ be a strongly semistable
vector bundle over a smooth curve $C$. Then there is a sequence of
genuinely ramified maps $f_i:D_i\to C$ such that $$Lim~
\frac{e_s(f_i^*(V))}{{\rm deg}~f_i} ~~=~~ \mu (V) $$ follows from
Lemma~\ref{genramofCIcurves}.  \end{remark}

Let $V_{\bullet}:=\{0=V_0\subset \cdots \subset V_l=V\}$ be the Harder
Narasimhan filtration of $V$. Let $d_i, r_i$ denote the degree and
rank of $V_i$ for $i = 0, 1, 2, \cdots l$.  Assume that the Harder
Narasimhan factors $V_i/V_{i-1}$ are strongly semistable if the
characteristic is positive.  By [La], any vector bundle on a curve in
positive characteristic has such a Harder Narasimhan Filtration after a
finite number of Frobenius pull backs. Now we can determine the asymptotic
slopes of $V$ for each $s$. 

\begin{theorem} \label{theorem2} $s\nu_s = d_i + (s-r_i)\mu_{i+1} $
where $r_i< s \leq r_{i+1}$ \end{theorem}

\begin{proof} First we note that if $L$ is any line bundle, then
$e_s(V\otimes L)=e_s(V)+ {\rm deg}~L$ and the spectrum ${\mathcal
{AS}}_s(V\otimes L)={\mathcal {AS}}_s(V) + {\rm deg} L$. Hence
$\nu_s(V\otimes L)= \nu_s(V) + {\rm deg} L$.  Since the
Harder-Narasimhan filtration of $V\otimes L$ is given by $V_i \otimes
L$, $d_i(V\otimes L) = d_i(V) + r_i\cdot {\rm deg} L$. Now assuming
the formula for $V\otimes L$, we get:

$$s(\nu_s(V\otimes L)) = d_i(V\otimes L) + (s-r_i) \mu_{i+1}(V\otimes L) $$
$$ s(\nu_s(V)) + s({\rm deg}~L) = d_i(V) + r_i({\rm deg}~L) + 
(s-r_i)(\mu_{i+1}(V)) + (s-r_i)({\rm deg}~L) $$
By simplifying, we obtain the formula for $V$. 

Hence by taking deg $L$ to be sufficiently
large we may assume that all $\mu_{i}$'s are positive.

Let $f:D\to C$ be any finite map and $W\subset f^*(V)$ be
a subbundle of rank $s$.  Let $W_j\subset f^*(V_j/V_{j-1})$ be the
saturation of the image of $W\cap f^*V_j$ in $f^*(V_j/V_{j-1})$.  Let
$s_j$ be the rank of $W_j$ and $\delta_j$ be equal to $\frac{{\rm
deg}~W_j}{{\rm deg}~f}$.  Then we have $s_j\leq r_j-r_{j-1}$, $r_i< \sum
s_j=s\leq r_{i+1}$ and $\frac{{\rm degree}~W}{{\rm deg}~f}\leq
\sum\delta_j$. Since $V_j/V_{j-1}$ is strongly semistable, we also have
$\delta_j\leq s_j\mu_j\leq (r_j-r_{j-1})\mu_j$ for all $j\geq 1$.
Now by comparing the first $i$ terms and the
rest of the terms we get the inequality:

$$\frac{{\rm deg}~W}{{\rm deg}~f}\leq \sum_{j=1}^l s_j\mu_j = \sum_{j=1}^i
s_j\mu_j+\sum_{j=i+1}^l s_j\mu_j \leq \sum_{j=1}^i (r_j-r_{j-1})\mu_j
+ \sum_{j=i+1}^l s_j\mu_{i+1} \leq d_i+ (s-r_i)\mu_{i+1} $$

To show the equality we produce subbundles of rank $s$ in coverings
with degree divided by the degree of the covering arbitrarily close to
$ d_i + (s-r_i)\mu_{i+1} $ where $r_i< s \leq r_{i+1}$. By
Theorem~\ref{theorem1}, we can find a covering $f:D\to C$ and a
subbundle $W_i\subset f^*(V_i/V_{i-1})$ of rank $s-r_i$ such that
$\mu(V_i/V_{i-1})-\mu(W_i)< \epsilon$. Let $W:=\pi^{-1}(W_i)$ be the
inverse image of $W_i$ in $V_i\subset V$ by the projection $\pi:V_i\to
V_i/V_{i-1}$. Then it shows that $d_i+ (s-r_i)\mu_{i+1}-\mu(W)/{\rm
deg~f}<\epsilon$. Hence the theorem.
\end{proof}

\begin{theorem}\label{theorem3} Let $f:Y\to X$ be a morphism and
$W\subset f^*V$ be a subbundle of rank $s$ with $r_i<s\leq r_{i+1}$,
such that $$\frac{{\rm deg}~(W)}{{\rm deg}~f}>d_i+(s-r_i)
\mu_{i+1}-(\mu(V_i)-\mu_{i+1})$$ Then $f^*V_i\subset W\subset
f^*V_{i+1}$.  
\end{theorem}

\begin{proof}
From the proof of Theorem~\ref{theorem2}, $s\nu_s=d_i+ (s-r_i)\mu_{i+1} $,
the inequalities becomes equality if and only if $s_j=r_j-r_{j-1}$ for all
$j$,
and hence $W\cap V_i=V_i$, proving $V_i\subset W$. 
\end{proof}

\section{The geometry of the spectrum}

Notice that the $s$-spectrum ${\mathcal {AS}}_s(V)\subset
[e_s(V),\overline\mu_{\rm max}(V)] \subset {\mathbb R}$ is a subset of the
bounded interval. Hence it has maximum, minimum, and cluster points.
We have described the supremum (asymptotic slopes) of the spectrum for
each $s$. This leads to the following natural question.

\noindent {\bf Question}: Are asymptotic $s$-slopes the only limit
points of the spectrum? Or is it likely to be dense in the interval
$[e_s(V),\mu_{\rm max}(V)]$?

Now we give a criterion for the asymptotic slopes to be an isolated
value for strongly semistable vector bundles.  

\begin{lemma} Let $V$ be a strongly semistable vector bundle. Then $\mu(V)$ is
an isolated point of the asymptotic $s$-spectrum ${\mathcal {AS}}_s(V)$ if
and only if ${\mathcal {AS}}_s(V)=\{~\mu(V)~\}$, i.e., there exists a
subbundle $W\subset V$ of rank $s$ such that $\mu(W)=\mu(V)$.
\end{lemma}

\begin{proof} Assume $W\subset V$ is a subbundle of slope
$\mu(V)$. Then for any map $f:D\to C$, $f^*(W)\subset f^*(V)$ is a
maximal subbundle and hence $\frac{\mu(f^*W)}{{\rm deg}~f}=\mu
(W)=\mu(V)$.  Hence the spectrum is a singleton, proving ${\mathcal
{AS}}_s(V)=\{~\mu(V)~\}$.

Let $V$ be a strongly semistable vector bundle such that $\mu(V)$ is
isolated in the spectrum ${\mathcal {AS}}_s(V)$.

From Remark~\ref{limitslope} and the hypothesis that $\mu(V)$ is
isolated, it follows that there is a genuinely ramified map such that
the pull back of $V$ has a subbundle of same slope. By
Lemma~\ref{stabletostable}~(c), this subbundle descends to a subbundle
of same slope as $V$.
\end{proof}

\end{document}